\documentclass[11pt,a4paper,final,reqno]{amsproc}
\usepackage{amssymb,epic,eepic,hyperref}

\providecommand{\BBb}[1]{{\mathbb{#1}}}

\providecommand{\cal}[1]{{\mathcal{#1}}}   
\newcommand{\goth}[1]{{\mathfrak #1}}

\newcommand{\Bcirc}{\overset{\lower 1.5pt%
              \hbox{$@,@,@,@,@,\scriptscriptstyle\circ$}}B{}}
\newcommand{\Binfty}{\overset{\lower 1.5pt%
              \hbox{$@,@,@,@,@,\scriptscriptstyle\infty$}}B{}}
\newcommand{\bigdot}{\mathbin{\raise.65\jot\hbox{$\scriptscriptstyle\bullet$}}}

\newcommand{\class}{\operatorname{class}}

\newcommand{\Dm}{\BBb D}

\newcommand{\dv}{\operatorname{div}}
\newcommand{\erd}{\overset{\lower 1pt\hbox{\large.}}{e}
                  \overset{\lower 1pt\hbox{\large.}}{r}}

\newcommand{\Fcirc}{\overset{\lower 1.5pt%
               \hbox{$@,@,@,@,@,\scriptscriptstyle\circ$}}F{}}
\newcommand{\fracc}[2]{{
                \textstyle\frac{#1}{\raise 1pt\hbox{$\scriptstyle #2$}}}}
\newcommand{\fracp}{\fracc1p}
\newcommand{\fracnp}{\fracc np}
\newcommand{\fracr}{\fracc1r}
\newcommand{\fracnr}{\fracc nr}
\newcommand{\fracci}[2]{{\frac{#1}{\raise 1pt\hbox{$\scriptscriptstyle #2$}}}}
\newcommand{\fracpi}{\fracci1p}
\newcommand{\fracri}{\fracci1r}

\newcommand{\grad}{\operatorname{grad}}

\newcommand{\mlap}{-\!\operatorname{\Delta}}

\newcommand{\op}[1]{\operatorname{#1}}

\newcommand{\N}{\BBb N}

\newcommand{\R}{{\BBb R}}
\newcommand{\Rn}{{\BBb R}^{n}}

\providecommand{\rom}[1]{\upn{#1}}

{
\end{minipage}%
\par\medskip\noindent}%
\newcounter{enmcount}\renewcommand{\theenmcount}{{\rm\arabic{enmcount}}}
\renewenvironment{enumerate}{%
\begin{list}{{\llap{\rm(\theenmcount)}}}{\setlength{\labelwidth}{\leftmargin}%
\usecounter{enmcount}}}{\end{list}}
\newcounter{rmcount}\renewcommand{\thermcount}{{\rm\roman{rmcount}}}

\newcounter{Rmcount}\renewcommand{\theRmcount}{{\rm\Roman{Rmcount}}}

\newcommand{\Z}{\BBb Z}

\renewcommand{\hat}[1]{\overset{{\scriptscriptstyle \wedge}}{#1}}

\numberwithin{equation}{section}

\newtheorem{thm}{Theorem}
\numberwithin{thm}{section}

\theoremstyle{definition}
\newtheorem{defn}[thm]{Definition}

 \numberwithin{exercise}{section}
 
\theoremstyle{remark}
\newtheorem{rem}[thm]{Remark}

\title[]{ ~ \\ ~ \\ ~ \\ ~ \\ ~ \\ ~ \\ ~ \\Regularity 
properties of semilinear 
boundary problems in Besov and Triebel-Lizorkin spaces}
\author[]{~ \\ ~ \\Jon Johnsen\\Mathematisches Institut\\
Friedrich--Schiller--Universit\"at Jena}
\address{Mathematical Institute, University of Copenhagen,
 Universitetsparken~5, DK-2100 Copenhagen~O, Denmark}
\email{jjohnsen@math.ku.dk}
\thanks{Supported by the Danish Natural Sciences Research
Counsil, no.~11--1221--1.
\\[6\jot]
{\tt Appeared in Journ\'ees ``\'Equations deriv\'ees partielles'', St.~Jean de Monts 1995 (Palaiseau, France 1995), pp.\ XIV1--XIV14.}}

\begin{document}

\maketitle
\vfill

\textbf{Abstract:} Semi-linear elliptic boundary problems with
non-linearities of product type are considered, in particular the
stationary Navier--Stokes equations. Regularity and existence results
are dealt with in the Besov and Triebel--Lizorkin spaces, and it is
explained how difficulties occurring for boundary conditions of a high
class may be handled.
\thispagestyle{empty}

\clearpage
\setcounter{page}{1}
\section{Summary}\label{summ-sect}
For simplicity's sake the following two model problems are considered 
on a bounded open set $\Omega\subset\Rn$, where $n\ge2$ and
$\Gamma:=\partial\Omega$ is $C^\infty$-smooth: first there is the 
Dirichl\'et problem
\begin{equation}
  \begin{aligned}
    \mlap u+u\partial_{x_1}u &= f &&\quad\text{in } &&\Omega,
    \\
    \gamma_0u &= \varphi &&\quad\text{on } &&\Gamma.
  \end{aligned}   
  \label{Dir-pb}
\end{equation}
Here $\gamma_0u=u|_{\Gamma}$ and 
$\mlap u=-(\partial^2_{x_1}+\dots+\partial^2_{x_n})u$.
Secondly there is the corresponding Neumann problem 
\begin{equation}
  \begin{aligned}
    \mlap u+u\partial_{x_1}u &= f &&\quad\text{in } &&\Omega,
    \\
    \gamma_1u &= \varphi &&\quad\text{on } &&\Gamma,
  \end{aligned}   
  \label{Neu-pb}
\end{equation}
where $\gamma_1u=\gamma_0(\vec n\cdot\grad u)$ with $\vec n$ denoting
the unit outward normal vectorfield near $\Gamma$. For the stationary
Navier--Stokes equations and other problems,
see Theorem~\ref{NS-thm} and Section~\ref{rem-sect}. 

The regularity of the solution $u(x)$ is studied here together with
the question of carrying over weak solutions to other spaces. To
obtain a unified treatment of various well-known scales of function
spaces, the Besov spaces $B^s_{p,q}$ are considered together with the
Triebel--Lizorkin spaces $F^s_{p,q}$; hereby $s\in\R$ and $p$ and
$q\in\,]0,\infty]$ in general, although $p<\infty$ is required
throughout for the $F^s_{p,q}$ spaces.

Among the various identifications, recall eg that
$B^{s}_{\infty,\infty}=C^s_*$ for $s>0$ (the H\"older--Zygmund spaces);
$B^{s}_{p,p}=W^s_p$ for $s\in\R_+\!\setminus\N$, $1<p<\infty$
(Sobolev--Slobodetskii);
$F^{s}_{p,2}=H^s_p$ for $s\in\R$,  $1<p<\infty$ 
    (Bessel--potentials) so in particular this encompasses the $W^k_p$ 
and  $L_p$;
$F^{0}_{p,2}=h_p$ for $0<p<\infty$ (local Hardy space).
The scales coincide when $p=q$, so $B^{s}_{2,2}=F^{s}_{2,2}=H^s$ is
the usual Sobolev space for $s\in\R$.

On $\Rn$ the spaces are defined by
means of Littlewood--Paley decompositions,
$B^s_{p,q}(\overline{\Omega})=r_\Omega B^s_{p,q}(\Rn)$ etc denotes 
the restriction to $\Omega$; on $\Gamma$ local coordinates are used.
A concise review of the definition and the properties of the Besov and
Triebel--Lizorkin spaces is given in \cite{JJ96ell}, so details are
omitted here; for a proper exposition the reader
is referred to the books of H.~Triebel \cite{T2,T3} and to
Theorems~3.6 and 3.7 in M.~Yamazaki's article \cite{Y1}.

\bigskip

For the Dirichl\'et problem above there is the following result:
\begin{thm} \label{Dir-thm}
  Let $u(x)$ in $F^s_{p,q}(\overline{\Omega})$ be a solution of
  \eqref{Dir-pb} for data $f(x)$ in $F^{t-2}_{r,o}(\overline{\Omega})$
  and $\varphi(x)$ in $B^{t-\fracri}_{r,r}(\Gamma)$, and suppose that 
  \begin{align}
    s&>\max(\tfrac{1}{2},\fracnp-1+\tfrac12\delta_{n2}),
    \label{Dir-cnd} \\
    t&>\max(\tfrac{1}{2},\fracnr-1+\tfrac12\delta_{n2}).
    \label{Dir-cnd'}  
  \end{align}
Then $u(x)$ is also an element of $F^{t}_{r,o}(\overline{\Omega})$.

Analogously, if $u\in B^s_{p,q}(\overline{\Omega})$, $f\in
B^{t-2}_{r,o}(\overline{\Omega})$ and $\varphi\in
B^{t-\fracri}_{r,o}(\Gamma)$, then \eqref{Dir-cnd}--\eqref{Dir-cnd'}
imply that $u\in B^t_{r,o}(\overline{\Omega})$.
\end{thm}

The conditions \eqref{Dir-cnd}--\eqref{Dir-cnd'} in the theorem are
natural, for both $\gamma_0$ and $B(v):=v\partial_1 v$ 
make sense on $B^s_{p,q}(\overline{\Omega})$ and 
$F^s_{p,q}(\overline{\Omega})$ when \eqref{Dir-cnd} holds. Actually
$B(\cdot)$ is even `better behaved' on these
spaces than $\mlap$ then; this is made precise below by
taking a specific $\delta=\delta(s,p)$ such that $\delta>0$.

If one denotes $\cal
A_D=\left(\begin{smallmatrix}\mlap\\ \gamma_0\end{smallmatrix}\right)$
and $\cal B(u)=\left(\begin{smallmatrix}B(u)\\
0\end{smallmatrix}\right)$, problem \eqref{Dir-pb} becomes 
$\cal A_Du+\cal B(u)=\left(\begin{smallmatrix}f\\ \varphi
\end{smallmatrix}\right)$. Then, if $\delta(s,p)>0$ and
$\delta(t,r)>0$, ie if $\cal B(\cdot)$ respects the {\em direct\/}
regularity properties of $\cal A_D$ at $(s,p,q)$ and $(t,r,o)$, the
theorem asserts that $\cal B(\cdot)$ also respects the {\em inverse\/}
regularity properties of $\cal A_D$ at these two parameters. Moreover,
this holds for both of the $B^s_{p,q}$ and $F^s_{p,q}$ scales.

For the Neumann problem $\cal A_Nu+\cal
B(u)=\left(\begin{smallmatrix}f\\ \varphi\end{smallmatrix}\right)$,
where $\cal A_N=\left(\begin{smallmatrix}\mlap\\ \gamma_1
\end{smallmatrix}\right)$, there~is 

\begin{thm} \label{Neu-thm}
  Let $u(x)$ in $F^s_{p,q}(\overline{\Omega})$ be a solution of
  \eqref{Neu-pb} for data $f\in F^{t-2}_{r,o}(\overline{\Omega})$
  and $\varphi\in B^{t-1-\fracri}_{r,r}(\Gamma)$, and suppose that 
  \begin{equation}
    s>\max(\fracp+1,\fracnp-1+\tfrac12\delta_{n2}),\quad
    t>\max(\fracr+1,\fracnr-1+\tfrac12\delta_{n2}).
    \label{Neu-cnd}
  \end{equation}
  Then $u(x)$ belongs to $F^{t}_{r,o}(\overline{\Omega})$.
  The analogous result holds in the $B^s_{p,q}(\overline{\Omega})$ spaces.
\end{thm}

It turns out that
Theorem~\ref{Neu-thm} is rather more complicated to prove than
Theorem~\ref{Dir-thm}. The reason for this is that the requirement
$s>\frac12$ is replaced by $s>\fracp+1$ (because of $\gamma_1$), 
which is `bad' in its $p$-dependence. Roughly speaking, this means that if
$F^s_{p,q}(\overline{\Omega})+F^t_{r,o}(\overline{\Omega})\subset
F^{s_1}_{p_1,q_1}(\overline{\Omega})$, then $(s_1,p_1,q_1)$ need not
satisfy \eqref{Neu-cnd} even if both $(s,p,q)$ and $(t,r,o)$ do so. 
As outlined in Section~\ref{Neu-sect} below the fine theory of
pointwise multiplication provides estimates of $B(\cdot)$, that may
be used to overcome the difficulties.

\bigskip

Instead of the model problems above, the methods may be applied to eg
the stationary Navier--Stokes equations. For each of the five boundary
conditions considered in \cite{GS2} one finds regularity results for
the solutions that correspond to either Theorem~\ref{Dir-thm} or
Theorem~\ref{Neu-thm}. See \cite[Thm.~5.5.5]{JJ93} or \cite{JJ95reg}
for this. 

In addition the existence of weak solutions
of the Dirichl\'et problem may be carried over to the $B^s_{p,q}$ and
$F^s_{p,q}$ spaces in this way. In more details the problem is:
\begin{equation}
  \begin{aligned}
    \mlap u+\smash{\sum_{j=1}^n} u\partial_j u +\grad\goth p&= f
    &\quad&\text{in } &&\Omega,
    \\
    \dv u&= g &\quad&\text{in } &&\Omega,
    \\
    \gamma_0 u&=\varphi &\quad& \text{on } &&\Gamma.
  \end{aligned}
  \label{NS-pb}
\end{equation}
Here the solution $(u,\goth p)$ and the data $(f,g,\varphi)$ are
sought such that 
\begin{equation}
  \begin{gathered}
   u\in B^{s}_{p,q}(\overline{\Omega})^n,\qquad 
   \goth p\in B^{s-1}_{p,q}(\overline{\Omega})
   \\
   f\in B^{s-2}_{p,q}(\overline{\Omega})^n,\qquad
   g\in B^{s-1}_{p,q}(\overline{\Omega}),\qquad
   \varphi\in B^{s-\fracpi}_{p,q}(\Gamma)^n,
  \end{gathered}
  \label{NS-spcs}
\end{equation}
for $s>\max(\frac 12,\fracnp-1+\frac12\delta_{n2})$; observe that the 
problem in \eqref{Dir-pb} may serve as a model problem for
\eqref{NS-pb}. For the $F^s_{p,q}$ spaces the requirement is the same,
but again $B^{s-\fracpi}_{p,q}(\Gamma)$ should be replaced by
$B^{s-\fracpi}_{p,p}(\Gamma)$. 

Concerning the existence of solutions when $g=0$ there is:

\begin{thm} \label{NS-thm}
  Let $\Omega\subset\Rn$, where $n=2$ or $3$, be a $C^\infty$-smooth
  open bounded set,
  and let $\Omega$ be connected with finitely many components of
  $\Gamma$, ie $\Gamma=\Gamma_1\cup\dots\cup\Gamma_N$.
 
  Suppose that the data $(f,0,\varphi)$ belong to the spaces indicated
  in \eqref{NS-spcs} for a parameter $(s,p,q)$ satisfying one of the
  following conditions:
 \begin{enumerate}
   \item[(1)] $s>\max(1,\fracc np+1-\frac n2)$;
   \item[(2)] $s>1$, $s=\fracc np+1-\frac n2$ and $q\le2$;
   \item[(3)] $s=1$ and $p\ge2\ge q$.
 \end{enumerate}
 Assume in addition that $\int_{\Gamma_j}\vec n\cdot\varphi=0$ for
 $j=1$,\dots,$N$. 

 Then there exists a solution $(u,\goth p)$ of \eqref{NS-pb} as in
 \eqref{NS-spcs} above.

 For the $F^s_{p,q}$ spaces the analogous result holds (for any
 $q\in\,]0,\infty]$ in \rom{(2)}).
\end{thm}
The special case with $(s,p,q)=(1,2,2)$ is identical to the classical
result on weak solutions, cf~\cite{Tm}. As a particular case
the theorem gives a solvability theory in the H\"older--Zygmund spaces
$C^s_*(\overline{\Omega})$ for $s>1$.

In addition solutions may be constructed by successive
approximations for any $s>\max(\frac12,\fracnp-1+\frac12\delta_{n2})$
in \eqref{NS-spcs} provided only that $\int_\Omega g=\int_\Gamma \vec
n\cdot\varphi$ and that the norms of the data are small enough; for
the present spaces, this is elaborated in \cite{JJ93}.  In comparison
Theorem~\ref{NS-thm} asserts that when $g=0$ and $s$ is sufficiently
large (plus some stricter conditions on $\Omega$ and $\varphi$), then 
solutions exist for arbitrarily large data. 

In view of this even the $C^s_*$ result should be new.

\bigskip

The purpose of this paper is only to indicate the proofs of the theorems;
a detailed exposition is in preparation \cite{JJ95reg}. The results are
based on \cite{JJ96ell, JJ94mlt}.

\section{The pseudo-differential boundary operators} \label{grn-sect}
For an efficient treatment of the problems in \eqref{Dir-pb}, 
\eqref{Neu-pb} and \eqref{NS-pb} one can utilise the
calculus of pseudo-differential boundary operators of L.~Boutet
de~Monvel \cite{BM71} for the linear parts. An extension of this 
calculus to the $B^s_{p,q}$ and
$F^s_{p,q}$ scales may be found in \cite{JJ96ell,JJ93} (with the
results of J.~Franke (partially contained) in \cite{F3} and the
$H^s_p$ and $B^s_{p,p}$ versions of G.~Grubb \cite{G3} as
forerunners).

Introductions to the calculus may be found in \cite[Sect.~2]{G2} and
\cite[Sect.~1.1 ff]{G1}, or \cite[Sect.~4]{G3}, so here it is recalled
that the generic object to study is a Green operator
\begin{equation}
   \cal A= \begin{pmatrix}P_\Omega+G&K \\[2\jot] T&S\end{pmatrix} \colon
  \begin{array}{ccc}
  C^\infty(\overline{\Omega})^N\\ \oplus\\ C^\infty(\Gamma)^M
  \end{array}
 \to
  \begin{array}{ccc}
  C^\infty(\overline{\Omega})^{N'}\\ \oplus \\ C^\infty(\Gamma)^{M'}
  \end{array},
 \label{grn1}  
\end{equation}
whereby $P_\Omega=r_\Omega Pe_\Omega$ denotes the truncation to
$\Omega$ of a pseudo-differential operator on $\Rn$; $T$ is a trace
operator, $K$ a Poisson operator and $S$ is a
pseudo-differential operator on $\Gamma$; finally $G$ is a singular
Green operator.

To assure that $P_\Omega(C^\infty(\overline{\Omega})^N)\subset
C^\infty(\overline{\Omega})^{N'}$ the so-called transmission condition
at $\Gamma$ is imposed on $P$ (cf~the elementary exposition in
\cite[Sect.~1]{G2}). More precisely, the results in 
\cite{JJ96ell} have been established for the space-uni\-formly estimated
calculus, for which the H\"ormander class $S^d_{1,0}(\Rn\times\Rn)$ is
the basic symbol class on $\Rn$; this version of the calculus has been
introduced systematically in \cite{GK}. Hence $P$ is required to
satisfy the uniform two-sided transmission condition at $\Gamma$, and
for $\cal A$ of the described kind the main result in \cite{JJ96ell}
is:

\begin{thm} \label{grn-thm}
  Suppose all entries in $\cal A$ have order $d\in\Z$ and that
  both $T$ and $P_\Omega+G$ are of class $r\in\Z$. Then there is
  continuity of
  \begin{align}
     \cal A &\colon
     \begin{array}{ccc}
       B^{s}_{p,q}(\overline{\Omega})^N\\ \oplus
       \\ B^{s-\fracpi}_{p,q}(\Gamma)^M 
     \end{array}
     \to
     \begin{array}{ccc}
       B^{s-d}_{p,q}(\overline{\Omega})^{N'}\\ \oplus \\   
       B^{s-d-\fracpi}_{p,q}(\Gamma)^{M'}
     \end{array},
     \label{grn2}  
  \\[2\jot]
     \cal A &\colon
     \begin{array}{ccc}
       F^{s}_{p,q}(\overline{\Omega})^N\\ \oplus
       \\ B^{s-\fracpi}_{p,p}(\Gamma)^M 
     \end{array}
     \to
     \begin{array}{ccc}
       F^{s-d}_{p,q}(\overline{\Omega})^{N'}\\ \oplus \\   
       B^{s-d-\fracpi}_{p,p}(\Gamma)^{M'}
     \end{array},
     \label{grn3}
  \end{align}
for $s>r+\max(\fracp-1,\fracnp-n)$. In both cases,
boundedness can only hold for $s<r+\max(\dots)$ if both $\class(T)$
and $\class(P_\Omega+G)$ are $<r$.
\end{thm}
When all symbols are poly-homogeneous and $\cal A$ is elliptic, the
theorem applies also to any parametrix $\widetilde{\cal A}$, and  it was
shown in \cite[Thm.~5.4]{G3} that $\widetilde{\cal A}$ can be taken of
class $r-d$. In general this result is best possible because $\cal A$
is a parametrix of $\widetilde{\cal A}$. (An exception is when $\cal A$
itself only contains a negligible part of class $r$.)

With obvious modifications the theorem also holds for multi-order and
multi-class operators (of the Douglis--Nirenberg type) or when either
$M$ or $M'=0$; see \cite[Thm.~5.2]{JJ96ell}. As examples there are then
$\cal A_D$ and $\cal A_N$; throughout 
\begin{equation}
  \widetilde{\cal A}_D:= \begin{pmatrix} R_D& K_D\end{pmatrix}=\cal
  A_D^{-1}
  \label{grn5}
\end{equation}
will serve as a special choice of parametrix of $\cal A_D$.

For convenience $\Dm_k$, with $k\in\Z$, will denote the admissible 
parameters $(s,p,q)$ for which the inequality
\begin{equation}
  s>k+\max(\fracp-1,\fracnp-n)
  \label{Dmk}
\end{equation}
holds. Equivalently this means that $s>k-1+\fracp+(n-1)(\fracp-1)_+$.

\section{Product estimates} \label{prd-sect}
The bilinear operator $B(v,w)=v\partial_1w$, that has been used above
with $B(v):=B(v,v)$, is analysed as the composite
\begin{equation}
  (v,w)\mapsto (v,\partial_1w)\mapsto \pi(v,\partial_1w),
\label{prd1}
\end{equation}
where $\pi(f,g)$ denotes $f(x)\cdot g(x)$. More precisely, 
$\pi(\cdot,\cdot)$ is the following
generalisation, that eg  allows $s>\frac12$ (instead of $s>1$) in
\eqref{Dir-cnd}:

\begin{defn} \label{pi-defn}
  For $u$ and $v\in \cal S(\Rn)$ let, with $\psi_k(\xi)=\psi(2^{-k}\xi)$,
  \begin{equation}
    \pi(u,v)=\lim_{k\to\infty} \cal F^{-1}(\psi_k\hat u)\cdot
    \cal F^{-1}(\psi_k\hat v),
    \label{prd2}
  \end{equation}
whenever the limit, calculated in $\cal D'(\Rn)$: (i) exists for each $\psi\in
C^\infty_0(\Rn)$ equal to $1$ near $0$, and (ii) is independent of such
$\psi$'s. 
\end{defn}

This product has been studied in \cite{JJ94mlt}, where it is shown
that it fills a part of the gap between two immediate meanings of
`pointwise multiplication': $\pi(f,u)=fu$ for $f\in\cal O_M$ and
$u\in\cal S'$, and $\pi(f_0,f_1)=f_0\cdot f_1$ when the $f_j$ lie in 
$L_{p_j}^{\op{loc}}\cap\cal S'$ such that
$0\le\fracc1{p_0}+\fracc1{p_1}\le1$ (so that $f_0\cdot f_1\in L_1^{\op{loc}}$).

Moreover, for an open set $\Omega\subset\Rn$ there is a restriction to
$\Omega$ defined as 
\begin{equation}
  \pi_\Omega(u,v)=\lim_{k\to\infty} r_\Omega(\cal F^{-1}(\psi_k\hat u_1)
  \cal F^{-1}(\psi_k\hat v_1)),
  \label{prd3}
\end{equation}
when the limit exists in $\cal D'(\Omega)$ and satisfies (i) and (ii) 
for some $u_1$ and $v_1$ in
$\cal S'$ such that $r_\Omega u_1=u$ and $r_\Omega v_1=v$. (The
existence of such a pair $(u_1,v_1)$ implies that the limit exists,
equals $\pi_\Omega(u,v)$ and fulfils (i) and (ii) for any other pair 
restricting to $(u,v)$.)

Perhaps more importantly, the continuity properties of
$\pi(\cdot,\cdot)$ may be obtained by para-multiplication. For spaces
over $\Omega$ the definition in \eqref{prd3} allows one
to carry boundedness over from $\pi(\cdot,\cdot)$ to
$\pi_\Omega(\cdot,\cdot)$, cf~\cite[Thm.~7.2]{JJ94mlt}.

For simplicity's sake only the needed $F^s_{p,q}$ results will be
recalled. For the Besov spaces it is necessary with a stricter
control over the sum-exponents $q$, but in the end this does not affect
the results in Theorems~\ref{Dir-thm}--\ref{NS-thm}; hence these
technicalities are omitted here.

\begin{thm} \label{prd-thm}
  The product in \eqref{prd3} is defined on
  $F^{s_0}_{p_0,q_0}(\overline{\Omega})\times
  F^{s_1}_{p_1,q_1}(\overline{\Omega})$ when 
  \begin{equation}
    s_0+s_1>\max(0,\fracc n{p_0}+\fracc n{p_1}-n),
    \label{prd4}
  \end{equation}
  and then there is boundedness 
  \begin{equation}
    \pi_\Omega (\cdot,\cdot)\colon
    F^{s_0}_{p_0,q_0}(\overline{\Omega})\oplus
    F^{s_1}_{p_1,q_1}(\overline{\Omega})
    \to F^{s_2}_{p_2,q_2}(\overline{\Omega})
    \label{prd5}
  \end{equation}
  if all of the following conditions are fulfilled:
  \begin{align}
    s_2&< \min(s_0,s_1);  
    \label{prd6}\\
    s_2-\fracc n{p_2}&\le\min(s_0-\fracc n{p_0},s_1-\fracc n{p_1},
          s_0+s_1-\fracc n{p_0}-\fracc n{p_1});  
          \label{prd7}\\
    s_2-\fracc n{p_2}&=s_0-\fracc n{p_0}\quad\text{and}\quad 
    s_1=\fracc n{p_1} \quad\text{hold only if}\quad p_1\le 1;
    \label{prd8}\\
    s_2-\fracc n{p_2}&=s_1-\fracc n{p_1}\quad\text{and}\quad 
    s_0=\fracc n{p_0} \quad\text{hold only if}\quad p_0\le 1.
    \label{prd9}
  \end{align}
Here it suffices with $s_2\le\min(s_0,s_1)$ in \eqref{prd6} provided
$q_2\ge q_j$ if $s_2=s_j$.
\end{thm}

For this result the reader is referred to the theorems in
\cite[Sect.s~6 and 7]{JJ94mlt}. Since $\pi_\Omega(\cdot,\cdot)$ is
commutative, it may be assumed that $s_0\ge s_1$, and then the value,
$p_1^*$, of $p_2$ for which there can be equality in both \eqref{prd6} and
\eqref{prd7} is given by the formula
\begin{equation}
  \fracc n{p^*_1}=\fracc n{p_1}+ (s_0-\fracc n{p_0})_-
  +(s_1-\fracc n{p_1}-(s_0-\fracc n{p_0})_+)_+.
  \label{prd10}
\end{equation}

\begin{rem}\label{vrtx-rem}
  In Theorem~\ref{prd-thm} the receiving spaces $F^{s_2}_{p_2,q_2}$ are
  determined implicitly by \eqref{prd6}--\eqref{prd9}. But, since
  $\Omega$ is bounded,
  $F^{s_1}_{p^*_1,q}(\overline{\Omega}) \hookrightarrow
  F^{s_2}_{p_2,q_2}(\overline{\Omega})$ holds in any case, if
  $q=q_1$ for $s_0>s_1$ and if $q=\max(q_0,q_1)$ for $s_0=s_1$.
  Thus the receiving space with $(s_2,p_2,q_2)=(s_1,p_1^*,q)$ may be
  considered as optimal.
\end{rem}

\section{The Dirichl\'et model problem} \label{Dir-sect}

This section concerns the proof of Theorem~\ref{Dir-thm}. Preference will
be given to the Triebel--Lizorkin spaces for simplicity, however, 
everything holds mutatis mutandem for the Besov spaces as well. 

\bigskip

Firstly, for the linear parts of \eqref{Dir-pb}, there is 
boundedness of 
\begin{equation}
  \cal A_D:=\begin{pmatrix} \mlap \\ \gamma_0\end{pmatrix}
  \colon F^s_{p,q}(\overline{\Omega})\to
  \begin{array}{c}  
   F^{s-2}_{p,q}(\overline{\Omega})\\ \oplus \\ F^{s-\fracpi}_{p,p}(\Gamma)
  \end{array}
  \label{i3}
\end{equation}
for each parameter $(s,p,q)$ with
$s>1+\max(\fracp-1,\fracnp-n)$, ie in $\Dm_1$.  

For $\cal A_D$, the calculus asserts that the parametrix 
$\widetilde{\cal A}_D$  is bounded in the opposite direction in \eqref{i3}
for each parameter $(s,p,q)\in\Dm_1$.

Secondly, when the non-linear term $u\partial_{x_1}u$ is taken into
consideration too, it is found
from \eqref{Dir-pb} that 
 \begin{equation}
   u=R_D f+ K_D\varphi- R_D(u\partial_{x_1}u).
 \label{i4}
 \end{equation}
This turns out to be meaningful when
$s>\max(\fracp,\tfrac12,\fracnp-\frac12)$ for $n=2$ and for
$s>\max(\fracp,\tfrac12,\fracnp-1)$ when $n\ge3$, so in general 
the condition is
\begin{equation}
  s>\max(\tfrac12,\fracnp-1+\tfrac12\delta_{n2}).
 \label{i5}
\end{equation}
To obtain this one can derive from Theorem~\ref{prd-thm} that
$u\mapsto (u,\partial_{x_1}u)\mapsto u\partial_{x_1}u$ is bounded,
for some $\delta(s,p)$,
\begin{equation}
  F^s_{p,q}(\overline{\Omega})\to F^s_{p,q}(\overline{\Omega})\oplus
  F^{s-1}_{p,q}(\overline{\Omega}) \to
  F^{s-2+\delta}_{p,q}(\overline{\Omega}),
 \label{i6}
\end{equation}
in general when $s>\max(\frac12,\fracnp-\frac{n-1}2)$,
cf~\eqref{prd4}. Here the deficit $\delta(s,p)$ measures how much the 
order of $B(\cdot)$ deviates from the order of $\mlap$.

It is essential that the non-linear term is {\em more\/} regular 
than $\cal A_D u$ when \eqref{i5} 
holds. In fact $\delta(s,p)$ equals $1+\min(0,s-\fracnp)$\,---\,an
increasing function of $s-\fracnp$\,---\,except that $\delta(\frac
np,p)=1-\varepsilon$ for some arbitrary $\varepsilon>0$. Hence
$\delta(s,p)>0$ as long as $s>\max(\frac12,\fracnp-1+\frac12\delta_{n2})$.

Thirdly, after these preparations, an iteration yields  that $u\in
F^t_{r,o}(\overline{\Omega})$: observe that in \eqref{i4} one has, by
\eqref{i6} and Theorem~\ref{grn-thm} applied to $\widetilde{\cal A}_D$, for the
summands on the right hand side that 
\begin{equation}
  R_Df+K_D\varphi\in F^t_{r,o}(\overline{\Omega}),\qquad
  R_D(u\partial_{x_1}u)\in F^{s+\delta(s,p)}_{p,q}(\overline{\Omega}).
 \label{i8}
\end{equation}
By determination of 
$F^{s_1}_{p_1,q_1}\supset F^t_{r,o}+F^{s+\delta(s,p)}_{p,q}$,
it follows that $u\in F^{s_1}_{p_1,q_1}$.
Application of \eqref{i6} then gives $R_D(u\partial_{x_1}u)\in
F^{s_1+\delta(s_1,p_1)}_{p_1,q_1}$ etc.

In the case $r=p$ one may take $p_1=p$, and, because $\delta(s,p)>0$, 
the process ends with the conclusion that $u\in F^{s}_{p,q}$ in
approximately $|t-s|/\delta(s,p)$ steps (as is well known).

For $r\ne p$ the conclusion follows  by consideration of four different
cases, namely those with the combinations of $s+\delta(s,p)\gtreqless t$ 
and $t-\frac nr\gtreqless s-\fracnp$.  The procedure is far easier 
to sketch with a diagram than with words, so the reader is 
referred to Figure~\ref{iter-fig}.

\begin{figure}[htbp]
\hfil
\setlength{\unitlength}{0.0125in}
\begin{picture}(401,245)(0,0)
\path(187,137)(130,80)\path(136,83)(130,80)(133,86)
\path(130,110)(100,80)\path(106,83)(100,80)(103,86)
\path(100,110)(70,80)\path(76,83)(70,80)(73,86)
\path(70,110)(55,95)
\path(195,140)(250,140)\path(244,142)(250,140)(244,138)
\path(250,160)(270,160)\path(264,162)(270,160)(264,158)
\path(270,180)(290,180)\path(284,182)(290,180)(284,178)
\path(290,200)(300,200)
\path(193.5,137.5)(210.5,120.5)\path(204.5,123.5)(210.5,120.5)(207.5,126.5)
\dottedline{5}(130,80)(130,110)
\dottedline{5}(100,80)(100,110)
\dottedline{5}(70,80)(70,110)
\dottedline{5}(250,140)(250,160)
\dottedline{5}(270,160)(270,180)
\dottedline{5}(290,180)(290,200)
\dottedline{5}(190,145)(190,165)
\dottedline{5}(190,175)(190,195)
\dottedline{5}(190,205)(190,215)
\dashline{5.000}(10,75)(205,75)(370,240)
\thicklines
\path(10,0)(10,240)\path(12.000,232.000)(10.000,240.000)(8.000,232.000)
\path(8,10)(360,10)\path(352.000,8.000)(360.000,10.000)(352.000,12.000)
\put(7,245){\makebox(0,0)[lb]{\raisebox{0pt}[0pt][0pt]%
{\shortstack[l]{{$s$}}}}}
\put(0,7){\makebox(0,0)[lb]{\raisebox{0pt}[0pt][0pt]%
{\shortstack[l]{{$0$}}}}}
\put(365,7){\makebox(0,0)[lb]{\raisebox{0pt}[0pt][0pt]%
{\shortstack[l]{{$\fracc{n}p$}}}}}
\put(186.5,139){\makebox(0,0)[lb]{\raisebox{0pt}[0pt][0pt]%
{\shortstack[l]{{$\scriptscriptstyle\times$}}}}}
\put(130,140){\makebox(0,0)[lb]{\raisebox{0pt}[0pt][0pt]%
{\shortstack[l]{{$(s+\delta,p)$}}}}}

\put(325,220){\makebox(0,0)[rb]{\raisebox{0pt}[0pt][0pt]%
{\shortstack[l]{{${(t,r)}$}}}}}
\put(330,220){\makebox(0,0)[cc]{\raisebox{0pt}[0pt][0pt]%
{\shortstack[l]{{${\scriptscriptstyle\times}$}}}}}
\put(120,220){\makebox(0,0)[lb]{\raisebox{0pt}[0pt][0pt]%
{\shortstack[l]{{$(t,r)\ \scriptscriptstyle\times$}}}}}
\put(30,60){\makebox(0,0)[lb]{\raisebox{0pt}[0pt][0pt]%
{\shortstack[l]{{$(t,r)$}}}}}
\put(30,78){\makebox(0,0)[lb]{\raisebox{0pt}[0pt][0pt]%
{\shortstack[l]{{$\scriptscriptstyle\times$}}}}}
\put(212,113.5){\makebox(0,0)[lb]{\raisebox{0pt}[0pt][0pt]%
{\shortstack[l]{{$\scriptscriptstyle\times$}}}}}
\put(202,100){\makebox(0,0)[lb]{\raisebox{0pt}[0pt][0pt]%
{\shortstack[l]{{$(t,r)$}}}}}
\put(187.6,168){\makebox(0,0)[lb]{\raisebox{0pt}[0pt][0pt]%
{\shortstack[l]{{$\scriptscriptstyle\circ$}}}}}
\put(187.6,198){\makebox(0,0)[lb]{\raisebox{0pt}[0pt][0pt]%
{\shortstack[l]{{$\scriptscriptstyle\circ$}}}}}
\put(205,60){\makebox(0,0)[lb]{\raisebox{0pt}[0pt][0pt]%
{\shortstack[l]{{$s=\max(\frac12,\fracnp-1)$}}}}}
\end{picture}
\caption{The $p$-dependent iteration (for $n=3$).}
\label{iter-fig}
\end{figure}

The figure displays the location of $F^{s+\delta(s,p)}_{p,q}$ and four
examples of $F^t_{r,o}$ corresponding to the
subdivision mentioned above. However, the sum-exponents $q$ and $o$
are not represented. The sector where $\delta(s,p)>0$, ie
$s>\max(\frac12,\fracnp-1+\tfrac12\delta_{n2})$, is indicated in 
dashed line; note that
the `$\times$' representing $F^{s+\delta}_{p,q}$ and $F^t_{r,o}$
all lie {\em inside\/} the sector.

The arrows in full
line indicates embeddings $F^{s+\delta(s,p)}_{p,q}\hookrightarrow
F^{s_1}_{p_1,q_1}$ etc: there are Sobolev embeddings down
the lines of slope $1$, and to the right along horizontal lines
because $\Omega$ has finite measure. The dotted line indicates the
improved knowledge of the non-linear term, ie the spaces
$F^{s_1+\delta(s_1,p_1)}_{p_1,q_1}$ etc. 

Note that one of the four cases is trivial since
$F^{s+\delta}_{p,q}\hookrightarrow F^t_{r,o}$, while another in
finitely many steps (indicated by `$\circ$') reduces to this or 
to one of the cases with either
the ``sawtooth'' or the ``staircase'' manoeuvres.

Altogether this leads to the proof of Theorem~\ref{Dir-thm}.

\begin{rem} 
  The procedure followed above has been used by S.~I.~Poho\-{\v z}aev,
  at least in the case with $s=t$ and $r>p$, cf \cite{Poh93}. 
\end{rem}

\section{The Neumann problem} \label{Neu-sect}
For the Neumann problem in \eqref{Neu-pb} the arguments in 
Section~\ref{Dir-sect} turn out to require more detailed estimates of
the non-linear term. The reasons for this will be described in the following.

\bigskip

For the problem in \eqref{Neu-pb}, one should take $(s,p,q)\in\Dm_2$, for then 
\begin{equation}
  \gamma_1\colon F^s_{p,q}(\overline{\Omega})\to F^{s-1-\fracpi}_{p,p}(\Gamma),
  \label{i11}
\end{equation}
is bounded. It is important here that $\gamma_1$ can not be
continuous from $F^s_{p,q}(\overline{\Omega})$ unless
$s\ge 2+\max(\fracp-1,\fracnp-n)$ holds, so the restriction for $(s,p,q)$
can not be essentially improved. For the Green operator 
$\cal A_N=\left(\begin{smallmatrix}\mlap\\
\gamma_1\end{smallmatrix}\right)$ this means that the class is $2$.

Since $\cal A_N$ is elliptic, the Boutet de~Monvel calculus asserts that
there exists a parametrix $\widetilde{\cal A}_N=\begin{pmatrix} R_N& K_N
\end{pmatrix}$ of class $0$\,---\,but not lower\,---\,that is bounded
\begin{equation}
  \widetilde{\cal A}_N\colon 
  \begin{array}{c}
  F^{s-2}_{p,q}(\overline{\Omega}) \\ \oplus \\ F^{s-1-\fracpi}_{p,p}(\Gamma) 
  \end{array}
  \to F^s_{p,q}(\overline{\Omega})
  \label{i12}
\end{equation}
when $(s,p,q)\in\Dm_2$. Hence the class of $R_N$ is $0$, so
$R_N$ can {\em not\/} be extended to an operator that is continuous from
$F^{s-2}_{p,q}(\overline{\Omega})$ when $s-2<\max(\fracp-1,\fracnp-n)$. 
Again the restriction on $(s,p,q)$ can not be essentially improved.

Contrary to the Dirichl\'et case above, $\widetilde{\cal A}_N$ is not an
inverse, but 
\begin{equation}
  \widetilde{\cal A}_N\cal A_N= 1-\cal R
  \label{i12'}
\end{equation}
for an operator $\cal R$ of order $-\infty$ and class $\le 2$. In
fact, this regularising operator may be taken as $\cal Ru=
|\Omega|^{-1}\int_\Omega u$ by a specific choice of $\widetilde{\cal A}_N$. 

Hence \eqref{i4} is
replaced by $u=R_N f+K_N\varphi+\cal R u-R_N(u\partial_{x_1}u)$, where
the first three terms belong to $F^t_{r,o}(\overline{\Omega})$. 

The conditions that make $u\partial_{x_1}u$ defined remain the same,
of course, so the assumption for $(s,p,q)$ in \eqref{i5} is here replaced by 
\begin{equation}
  s>\max(\fracp+1,\fracnp-1+\tfrac12\delta_{n2})
  \label{i13}
\end{equation}
and similarly for $(t,r,o)$. Cf~\eqref{Neu-cnd}.

Now the cases with $t<s+\delta(s,p)$ and $t-\frac nr>s+\delta-\frac
np$ are rather more complicated than the corresponding cases for the
Dirichl\'et problem. The reason for this is that the space
$F^{s_1}_{p_1,q_1}\supset F^{s+\delta}_{p,q}+F^t_{r,o}\ni u$ seemingly
may be much too large for an application of the non-linear operator $v\mapsto
R_N(v\partial_{x_1}v)$ to it. 

Indeed, whilst $t=s_1$ the integral-exponent $p_1$ is smaller than
$r$, and in many cases $s_1<1+\frac1{p_1}$ holds, 
although $t>1+\frac1r$. An example is sketched in Figure~\ref{*par-fig} 
below, where the sector determined by \eqref{i13} is indicated by
dashes. 
When $s_1$ is close to $\fracc n{p_1}-1$, then the deficit
$\delta(s_1,p_1)$ is close to $0$, and so  eventually
\begin{equation}
  s_1-2+\delta(s_1,p_1)<\fracc1{p_1}-1.
  \label{i14}
\end{equation}
In such cases, since $\class(R_N)=0$, the solution operator $R_N$ simply 
{\em does not\/} make sense on $F^{s_1-2+\delta(s_1,p_1)}_{p_1,q_1}
(\overline{\Omega})$, as recalled after \eqref{i12}. 
By comparison with the Dirichl\'et problem, the iteration is seemingly
unable to begin.

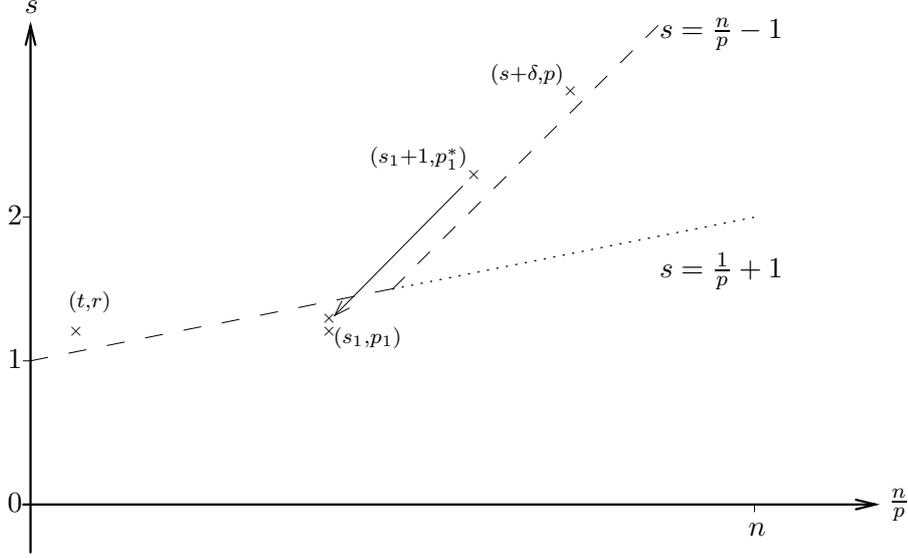
\begin{figure}[htbp]
\hfil
\setlength{\unitlength}{0.0125in}
\begin{picture}(401,246)(0,0)
\path(310,20)(310,17)
\path(10,80)(7,80)
\path(10,140)(7,140)
\dashline{8.000}(10,80)(160,110)(270,220)
\dottedline{4}(160,110)(310,140)
\path(189,153)(136,99)\path(142,102)(136,99)(139,105)
\thicklines
\path(10,0)(10,220)
\path(12.000,212.000)(10.000,220.000)(8.000,212.000)
\path(7,20)(360,20)
\path(352.000,18.000)(360.000,20.000)(352.000,22.000)
\put(0,17){\makebox(0,0)[lb]{$0$}}
\put(365,20){\makebox(0,0)[lc]{$\fracc{n}p$}}
\put(7,225){\makebox(0,0)[lb]{$s$}}
\put(307,7){\makebox(0,0)[lb]{$n$}}
\put(25,90){\makebox(0,0)[lb]{$\scriptscriptstyle\times$}}
\put(190,155){\makebox(0,0)[lb]{$\scriptscriptstyle\times$}}
\put(230,190){\makebox(0,0)[lb]{$\scriptscriptstyle\times$}}
\put(130,95){\makebox(0,0)[lb]{$\scriptscriptstyle\times$}}
\put(130,90){\makebox(0,0)[lb]{${\scriptscriptstyle\times}$}}
\put(135,85){\makebox(0,0)[lb]{$\scriptstyle (s_1,p_1)$}}
\put(150,160){\makebox(0,0)[lb]{$\scriptstyle (s_1+1,p^*_1)$}}
\put(200,195){\makebox(0,0)[lb]{$\scriptstyle (s+\delta,p)$}}
\put(25,100){\makebox(0,0)[lb]{$\scriptstyle (t,r)$}}
\put(0,77){\makebox(0,0)[lb]{$1$}}
\put(0,137){\makebox(0,0)[lb]{$2$}}
\put(270,210){\makebox(0,0)[lb]{$s=\fracc np-1$}}
\put(270,110){\makebox(0,0)[lb]{$s=\fracp+1$}}
\end{picture}
\hfil
\caption{The integral-exponent $p^*$.} 
 \label{*par-fig}
\end{figure}

At this place the fine theory of pointwise multiplication offers a
remedy. In fact, one can do better than regarding the non-linear
operator $v\mapsto v\partial_{x_1}v$ as one of order $2-\delta(s,p)$,
as in \eqref{i6} above. 

The problem only arises when $\delta(s,p)$ is close to $0$,
hence only for $s<\fracc np$, and then $v\partial_1 v$ may be seen to
factor through a space with smoothness index $\bold{s-1}$. More
exactly, with $\tfrac{n}{p^*}=\fracnp+(\fracnp-s)_+$, 
Theorem~\ref{prd-thm} gives that
\begin{equation}
  B(\cdot) \colon F^s_{p,q}(\overline{\Omega})
  \to  F^{s-1}_{p^*,q}(\overline{\Omega})
  \label{i16}
\end{equation}
is a bounded non-linear operator for $s>\max(\frac12,\fracc
np-1+\tfrac12\delta_{n2})$. With $(s,p,q)$ in the subsector given by 
$s>\max(1,\fracc np-1+\tfrac12\delta_{n2})$, in
which the $(s_1,p_1,q_1)$ above lies, it may be checked that 
the receiving space $F^{s-1}_{p^*,q}$ in \eqref{i16} lies
{\em above\/} the critical broken line $s=\max(\fracp-1,\fracnp-n)$,
so that $(s-1,p^*,q)\in\Dm_0$. See \cite{JJ95reg} or \cite{JJ93} for this.

According to \eqref{i12} this assures that $R_N$ may be applied to
$F^{s-1}_{p^*,q}(\overline{\Omega})$, so, because $s-1-\frac
n{p^*}=s+\delta(s,p)-\fracnp$, there is (after all) boundedness of  
\begin{equation}
  F^s_{p,q}(\overline{\Omega})\xrightarrow{B(\cdot)}
  F^{s-1}_{p^*,q}(\overline{\Omega}) \xrightarrow{R_N}
  F^{s+1}_{p^*,q}(\overline{\Omega}) \hookrightarrow 
  F^{s+\delta(s,p)}_{p,q}(\overline{\Omega})
  \label{i17}  
\end{equation}
for {\em all\/} $s>\max(1,\fracc np-1+\tfrac12\delta_{n2})$. 
Evidently this last sector is
stable under the forming of the intermediate parameters $(s_1,p_1,q_1)$,
$(s_2,p_2,q_2)$,\dots.

In principle also the cases with $t>s+\delta$ and
$t-\fracnr<s+\delta-\fracnp$ need a special argument, but also here
\eqref{i17} may be applied. Altogether the iteration used for 
the Dirichl\'et problem applies also to the Neumann problem.
\medskip

\section{Final remarks} \label{rem-sect}
(1) To prove Theorem~\ref{NS-thm}, notice that each of the conditions
(1)--(3) there implies that the data belong to the spaces
considered in Theorem~2.1 of \cite[App.~1]{Tm}. Hence there is a weak
solution to which the regularity results apply. For details, see
\cite[Thm.~5.5.5]{JJ93} or \cite{JJ95reg}.

(2) The iteration methods apply also to the von Karman equations for a
plate in $\Omega\subset\R^2$, or to problems with a suitable
semi-linear perturbation of an injectively elliptic Green operator
$\left(\begin{smallmatrix}P_\Omega+G& K\\ T&S\end{smallmatrix}\right)$
in the calculus.
\linebreak[5]


\end{document}